\documentclass{baustms}

\citesort

\theoremstyle{cupthm}
\newtheorem{thm}{Theorem}[section]

\newtheorem{lemma}[thm]{Lemma}
\theoremstyle{cupdefn}
\newtheorem{defn}[thm]{Definition}
\theoremstyle{cuprem}

\numberwithin{equation}{section}

\begin{document}
\runningtitle{A note on the hyperbolicity of the non-wandering sets of real quadratic maps}
\title{A note on the hyperbolicity of the non-wandering sets of real quadratic maps}
\author[1]{Diyath Pannipitiya}
\address[1]{LD 257, Department of Mathematical Sciences, Indiana University - Purdue University -- Indianapolis. IN. 46202. United States.\email{dinepann@iu.edu}}

\authorheadline{D. N. Pannipitiya}


\support{This paper is partially supported by NSF (Award Number: 2154414) }

\begin{abstract}
The goal of this paper is to discuss about the hyperbolicity of the non-wandering set $\mathcal{NW}(f_c)$ of real quadratic function $f_c(x)=x^2+c$ when $c\in (-\infty, -2]$. Even though the results we present here are not new, it is not easier to find the proofs of them. We present two different ways to prove the hyperbolicity of $\mathcal{NW}(f_c)$ for the``considerably difficult case'' \cite{bib2} of when $c$ is closer to $-2$.
\end{abstract}


\maketitle

\section{Intnroduction}
 
\begin{defn}
    (Non-wandering set)\\
 The non-wandering set of a function $f:X\longrightarrow X$ is defined as 
 \[\mathcal{NW}(f):= \{x: \text{for any nbhd $U$ of $x$, there is $n\in \mathbb{N}$ such that $f^{on}(U)\cap U \neq \emptyset$}\}.\]
\end{defn}
 Here, usually, $X$ is a compact metric space and $f$ is a continuous function.
 
 \begin{defn} (Hyperbolicity)\\
 Let $X$ be a metric space and let $f:X\longrightarrow X$ be a diffeomorphism. Let $M\subset X$ be an $f$-invariant set. We say $M$ is hyperbolic if there exist constants $c, C>0$ and $\lambda > 1$ for all $x\in M$  there is a decomposition $T_xM = E_x^s \bigoplus E_x^u$ of the tangent space of $M$ at $x$, such that $D_x(E^s_x) = E^s_{f(x)}$ and $D_x(E^u_x) = E^u_{f(x)}$. Here \[E_x^s:= \{v\in M: ||(D_xf^{on})v|| \leq c\Lambda ^{-n} ||v||\ \text{for all $n\in \mathbb{N}$}\}.\] and \[E_x^u:= \{v\in M: ||(D_xf^{on})v|| \geq C\Lambda^n  ||v||\ \text{for all $n\in \mathbb{N}$}\}\]
 are the stable and unstable manifollds of $x$ respectively.
   \end{defn}

\section{$\mathcal{NW}(f_{-2})$ is not hyperbolic.}
 
 First we show that $(f_{-2})|_{[-2,2]}$ is semi-conjugate to the shift map $\sigma$ on $\{0,1\}^{\mathbb{N}}$ - the one sided sequence space on two symbols. It is clear that ${f_{-2}([-2,2]) = [-2,2]}$. Let $[a,b]\subset [-2,2]$. Then 
 \begin{align}
     |f([a,b])|=\sqrt{b+2}-\sqrt{a+2}=\frac{b-a}{\sqrt{b+2}+\sqrt{a+2}}\leq \frac{|[a,b]|}{2\sqrt{2}} \label{eq1}.
 \end{align}
 Therefore by denoting the intervals $[-2,0]$ and $[0,2]$ by $0$ and $1$ respectively, we get that the map $\pi: \{0,1\}^{\mathbb{N}}\longrightarrow [-2,2] $ defined by $\pi (\theta)=x$, where $\theta$ is the $\{0,1\}$-name of $x$, is a semi-conjugacy of $\sigma$ with $f$ (that is $\pi \circ \sigma = f \circ \pi$). The reason is; the inequality \eqref{eq1} implies that each sequence in $\{0,1\}^{\mathbb{N}}$ is a $\{0,1\}$-name of only one point in $[-2,2]$. Thus $f$, when restricted to $[-2,2]$, is a factor of $\sigma$.\\
 
 Because there are points in $\{0,1\}^{\mathbb{N}}$ whose orbits are dense in $\{0,1\}^{\mathbb{N}}$ (under $\sigma$), the semi-conjugacy says that  we can find
 \begin{align}
     \text{$x_0\in [-2,2]$ such that $\mathcal{O}_{f_{-2}^+}(x_0)$ is dense in $[-2,2]$.}
 \end{align}
 
 Therefore for any interval $I\subset [-2,2]$, there is $m\in \mathbb{N}$ such that $f_{-2}^m(I)\cap I\neq \emptyset$. This says $\mathcal{NW}(f_{-2}) = [-2,2]$. 
 
 Now we show $\mathcal{NW}(f_{-2})$ is not hyperbolic. From the previous part, $0\in \mathcal{NW}(f_{-2})$. And it is clear that $D_0f_{-2}^n=0$ for any $n\in \mathbb{N}$. Thus $E_0^s=\mathbb{R}$ and $E_0^u=\{0\}$. Thus if $\mathcal{NW}(f_{-2})$ is hyperbolic, then we must have $D_0(E_0^s)=E^s_{f_{-2}(0)}=E^s_{-2}$. Notice that $f^2_{-2}(0)=f_{-2}(-2)=2$ is a fixed point. Hence $||D_{-2}f^n_{-2}||=2^{2n}$ for any $n\in \mathbb{N}$. This implies $E_{-2}^u=\mathbb{R}$, and therefore $E^s_{-2}=\{0\}$. Because $D_0(E_0^s)=\mathbb{R}$, we then have $D_0(E_0^s)\neq E^s_{-2}$. Contradiction! Therefore $\mathcal{NW}(f_{-2})$ is not hyperbolic.

\section{$\mathcal{NW}(f_c)$ is hyperbolic for all $c<-2$.}
Let $c<-2$. For simplicity, let $f=f_{c}$. Let $p_-$ and $p_+$ be the two fixed points of $f$ such that $p_-<0<p_+$. Let $\alpha\in [0,2]\backslash \{p_+\}$ be such that $f^2(\alpha)=p_+$ (see \eqref{fig1}). Then notice that the orbit of any point in the interval $(-\alpha, \alpha)$ escapes to infinity under the iterations of $f$. Let $J=[-p_+,p_+]\backslash (-\alpha, \alpha)$. Let $C=\bigcap_{n=0}^\infty f^{-n}(J)$. So, $\mathcal{NW}(f)=C$. \\

Notice that if $\alpha$ is large enough so that $|D_xf|>1$ for all $x\in [-p_+,p_+]\backslash (-\alpha, \alpha)$ (more precisely when $\alpha > \frac{1}{2}$), then for each $x\in C$ and for each $n\in \mathbb{N}$ we have that, \[|D_xf^n|=|(2x)^n|\geq (2\alpha)^n\geq \alpha 2^n.\] Thus $E^u_x=\mathbb{R}$, $E^s_x=\{0\}$, and $f|_C$ is a factor of the shift map $\sigma$ on $\{0,1\}^{\mathbb{N}}$. This says $\mathcal{NW}(f)$ is in-fact a Cantor set.  Therefore it is not hard to prove the hyperbolicity of $\mathcal{NW}(f)$ when $\alpha>1/2$ \Big(equivalently when $c < -\frac{5+2\sqrt{5}}{4}$\Big).\\

The hard part is to show the hyperbolicity when $0<\alpha \leq \frac{1}{2}$ \Big(equivalently when $-\frac{5+2\sqrt{5}}{4}\leq c < -2$\Big).\\

 Let $c\in [-\frac{5+2\sqrt{5}}{4},-2)$. Let $p=p_+$ be the positive fixed point of $f$. Then the end-points of $\mathcal{NW}(f)$ is $p$ and $-p$.
We show the hyperbolicity of $\mathcal{NW}(f)$ using two different methods.\\

\underline{\textbf{Method A}. (Using the hyperbolic metric)}\\

First let's recall some definitions.
\begin{defn}
Let $J\subset T$ be open and bounded intervals in $\mathbb{R}$ and let $L, R$ be the two intervals of $T\backslash J$. Define the cross-ratio of $T$ and $J$ as
\[D(T,J):= \frac{|J| |T|}{|L| |R|}.\]
\end{defn}

\begin{defn}
Let $J=(x,y)\subset T$ be open bounded intervals and let $L, R$ be the two intervals of $T\backslash J$. Define the hyperbolic metric on $T$ as 
\[\rho_T(x,y):= \log \frac{|L\cup J| |J \cup R|}{|L| |R|}.\]
\end{defn}

\textbf{Remark.} It is not hard to see that
\begin{align}
    \rho_T(x,y) = \log (1+ D(T,J) ).\label{3.1}
\end{align}

\begin{lemma}
If $J\Subset T$, then there exist $\Lambda \gneq 1$ for all $x,y\in J$ such that \[\rho_J(x,y)\geq \Lambda \rho_T(x,y).\]
\end{lemma}

\begin{proof}
Let $x,y\in J$. Let $M=(x,y)$. Let $P,Q$ be the the two disjoint intervals of $J\backslash (x,y)$. And let $L,R$ be the two disjoint intervals of $T\backslash J$.

\begin{center}
    \includegraphics[scale=0.4]{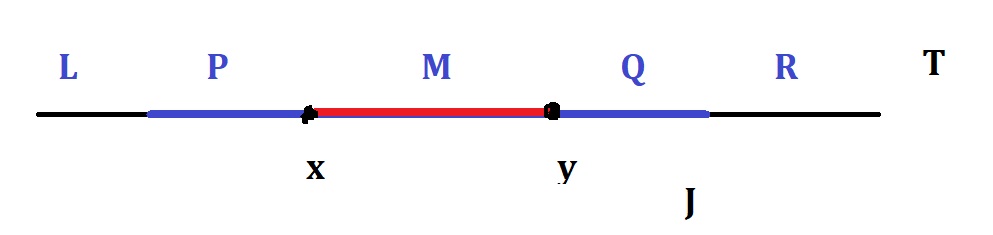}
\end{center}

Then
\[\frac{\rho_J(x,y)}{\rho_T(x,y)}=\frac{\ln\frac{(|P|+|M|)(|M|+|Q|)}{|P||Q|}}{\ln\frac{(|L|+|P|+|M|)(|M|+|Q|+|R|)}{(|L|+|P|)(|Q|+|R|)}}.\]

So this is strictly greater than 1 iff \[\frac{\frac{(|P|+|M|)(|M|+|Q|)}{|P||Q|}}{\frac{(|L|+|P|+|M|)(|M|+|Q|+|R|)}{(|L|+|P|)(|Q|+|R|)}}\gneq 1.\]
iff
\[\frac{|P|+|M|}{|L|+|P|+|M|}\cdot \frac{|M|+|Q|}{|M|+|Q|+|R|}\gneq \frac{|P|}{|L|+|P|}\cdot \frac{|Q|}{|Q|+|R|}.\]

iff
\[\frac{|L|+|P|}{|P|}\cdot \frac{|Q|+|R|}{|Q|}\gneq \frac{|L|+|P|+|M|}{|P|+|M|}\cdot \frac{|M|+|Q|+|R|}{|M|+|Q|}.\]

iff
\[\frac{|L|}{|P|}+\frac{|R|}{|Q|}+\frac{|L||R|}{|P||Q|}\gneq \frac{|L|}{|P|+|M|}+\frac{|R|}{|M|+|Q|}+\frac{|L||R|}{(|P|+|M|)(|M|+|Q|)}.\]

Which is true when $|M|\neq 0$. If $|M|=0$, that is if $x=y$, then by L'Hospital's rule we get
\[\frac{\rho_J(x,y)}{\rho_T(x,y)}=\lim_{|M|\longrightarrow 0}\frac{\ln\frac{(|P|+|M|)(|M|+|Q|)}{|P||Q|}}{\ln\frac{(|L|+|P|+|M|)(|M|+|Q|+|R|)}{(|L|+|P|)(|Q|+|R|)}} = \frac{|P|+|Q|}{|P||Q|}\cdot \frac{(|L|+|P|)(|Q|+|R|)}{|L|+|P|+|Q|+|R|}.\]
This is strictly greater than 1 iff
\[\frac{|R|}{|Q|}+\frac{|L|}{|P|}+\frac{|L||R|}{|P||Q|}\gneq \frac{|L|}{|P|+|Q|} + \frac{|R|}{|P|+|Q|}.\]

Which is true as $|L|^2+|R|^2\neq 0$.\\
Thus we have $\frac{\rho_J(x,y)}{\rho_T(x,y)}\gneq 1$. Because $J$ is compactly contained in $T$,  \[\Lambda:=\min_{x,y \in J} \frac{\rho_J(x,y)}{\rho_T(x,y)}\gneq 1 \]
exists. Hence we have the claim.
\end{proof}

\begin{defn}
Let $g:T \longrightarrow \mathbb{R}$ be a continuous and monotone function and let $J\subset T$. Define the cross-ratio of $g$ as
\[B(g,T,J):= \frac{D(g(T),g(J))}{D(T,J)}.\]
\end{defn}

\textbf{Remark.}
If $J$ and $T$ has common boundary points, then we take $\limsup{B(g,T,J_n)}$ where $(J_n)_{n=1}$ is an increasing sequence of nested sub-intervals of $J$ such that $J_n\nearrow J$.

\begin{defn}
For a $C^3$ map $g:T\subset \mathbb{R} \longrightarrow \mathbb{R}$, $\text{if}\ Dg(x)\neq 0$, define the Schwarzian derivative as
\[Sg(x):=\frac{D^3g(x)}{Dg(x)}-\frac{3}{2}\Big( \frac{D^2g(x)}{Dg(x)}\Big)^2.\]
\end{defn}

For the simplicity of this paper, we state some results from \textbf{De Melo and Van Strein}'s book \textbf{One Dimensional Dynamics} \cite{bib2} without proofs.\\

\textbf{Property 4} (Pg 273). If $g: T\longrightarrow \mathbb{R}$ is a $C^3$ map with $Sg<0$, then
\[B(g,T^*,J^*)\gneq 1\] for all pairs of intervals $J^*\subset T^* \subset T$.
Here we assume $g$ to be monotone on $T$.\\

Letting $T^*=T$ in this property, by \eqref{3.1}, we get

\begin{lemma}
Let $g: T\longrightarrow \mathbb{R}$ be a monotone function with $Sg<0$ and $Dg\neq 0$. Then \[\rho_{g(T)}(g(x), g(y))\geq \rho_T(x,y)\] for all $x,y\in T$.
\end{lemma}

From $Lemma\ 3.3$ and $Lemma\ 3.6$, there is $\Lambda \gneq 1$ for any $x,y\in J$ such that
\begin{align}
   \rho_J(f(x),f(y))\geq \Lambda \rho_J(x,y). \label{2}
\end{align}
\\

Now let's go back to $f(x)=x^2+c$ where $x\in \mathcal{NW}(f)$ (so here we restrict $f$ to its non-wandering set $\mathcal{NW}(f)$). Recall $0<\alpha\leq \frac{1}{2}$ is the positive second pre-image of $p$, where $p$ is the positive fixed point of $f$. Thus it is the minimum positive number in $\mathcal{NW}(f)$. Hence $\mathcal{NW}(f)$ is contained symmetrically in $[-p,-\alpha]\sqcup [\alpha, p]$. Let $I_L:=[-2p,\frac{-\alpha}{2}]$ and $I_R:=[\frac{\alpha}{2}, 2p]$. Thus $I_L\supsetneq \mathcal{NW}(f)\cap [-p, -\alpha]$ and $I_R\supsetneq \mathcal{NW}(f)\cap [\alpha, p]$.\\

Let $n\in \mathbb{N}$ be fixed and let $x\in \mathcal{NW}(f)$. WLOG suppose $x\in I_R$. Let $\Lambda \gneq 1$ be such that $\eqref{2}$ holds for $f^n$ on $I_R$. Let $h>0$ be sufficiently small such that $f^i(x-h)$ and $f^i(x+h)$ are in the same $I_L$ or $I_R$ as of $f^i(x)$, for all $i\in \{1,\cdots n\}$.
          \begin{figure}[h]\label{fig1}
        \centering
  \includegraphics[width=0.8\linewidth]{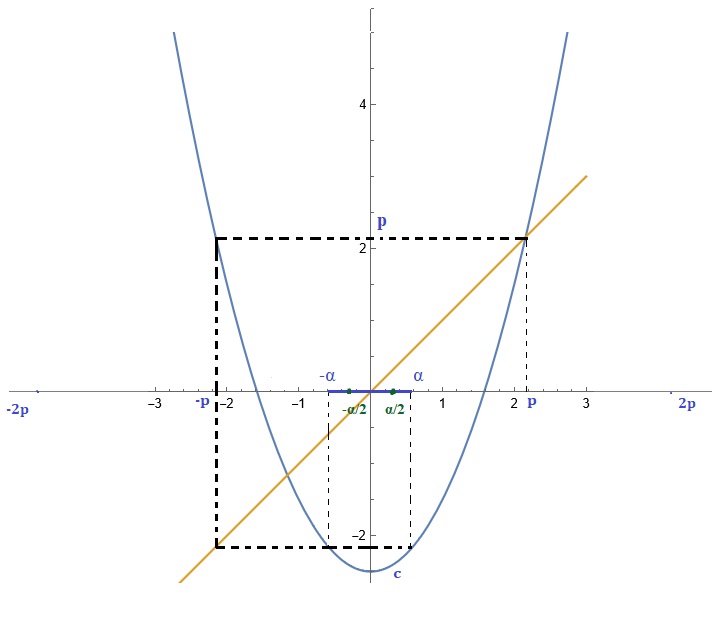}
        \caption{Graph of $f_c(x)=x^2+c$ when $c<-2$}
            \end{figure} 
Notice that on each $I_R, I_L$, the function $f$ is monotone. Also notice that because $f$ is an even function, $D(I_R, (f^i(x),fI(x+h))) = D(I_L, (f^i(x),fI(x+h)))$. And $d_{I_R}(x+h,x)=d_{I_L}(x+h,x)$. Thus by $(3)$, $d_{I_R}(f^n(x+h),f^n(x)) \geq \Lambda^n d_{I_R}(x+h,x)$. This implies $1+D(I_R, (f(x),f(x+h)))\gneq [1+D(I_R, (x,x+h)]^{\Lambda^n}$. \\
Which implies $D(I_R, (f(x),f(x+h)))\gneq \Lambda^n D(I_R, (x,x+h)$. Hence
\begin{align*}
    \Lambda^n &\lneq \frac{D(I_R, (f^n(x),f^n(x+h)))}{D(I_R, (x,x+h)}\\
    &=\frac{|f^n(x+h)-f^n(x)|(2p-\frac{\alpha}{2})}{|f^n(x)-\frac{\alpha}{2}||2p-f^n(x+h)|}\Big/ \frac{h(2p-\frac{\alpha}{2})}{(x-\frac{\alpha}{2})(2p-x-h)}\\
    &=\frac{|f^n(x+h)-f^n(x)|}{h}\cdot \frac{(x-\frac{\alpha}{2})(2p-x-h)}{|f^n(x)-\frac{\alpha}{2}||2p-f^n(x+h)|}.
    \end{align*}
Taking the limit as $h\longrightarrow 0$ (limit exists as $f\in C^\infty)$,
\begin{align*}
     \Lambda^n &\lneq |Df^n(x)| \frac{(x-\frac{\alpha}{2})(2p-x)}{|f^n(x)-\frac{\alpha}{2}||2p-f^n(x)|}.
\end{align*}
Because the end-points of $I_L$ and $I_R$ are not in $\mathcal{NW}(f)$, we then have
\[ \Lambda^n \lneq |Df^n(x)| \frac{(2p-\frac{\alpha}{2})^2}{(\alpha-\frac{\alpha}{2})(2p-p)}= |Df^n(x)| \frac{2(2p-\frac{\alpha}{2})^2}{\alpha p}.\]

Let $c=\frac{\alpha p}{2(2p-\frac{\alpha}{2})^2}$. Since $c$ does not depend on neither $x$ nor $n$, we get 
\[|Df^n(x)|> c\lambda^n,\ \text{for all}\ x\in \mathcal{NW}(f), n\in \mathbb{N}.\]

This proves $\mathcal{NW}(f)$ is hyperbolic.\\

\underline{\textbf{Method B}. (Using a weighted norm)}\\

Suppose we were able to find a continuous function $\omega: \mathcal{NW}(f) \longrightarrow (0,\infty)$ such that the weighted norm on $TM_x$ defined by \[||\overrightarrow{v}||_w:=\omega(x)\cdot ||\overrightarrow{v}||\] has the property: 
\begin{align}
    \text{there exists $\lambda \gneq 1$ for all $x\in \mathcal{NW}(f)$},\ ||Df(x)\overrightarrow{v}||_w \geq \Lambda ||\overrightarrow{v}||_w \text{for all $\overrightarrow{v}\in TM_x$.}
\end{align}

Because $\omega$ is non-vanishing on $\mathcal{NW}(f)$ and $\mathcal{NW}(f)$ is compact, we can find $\omega_0$ and $\omega_1$, say, such that $\omega_1\geq \omega(x) \geq \omega_0 \gneq 0$ for all $x\in \mathcal{NW}(f)$. Let $n\in \mathbb{N}$. Then
\begin{align*}
    ||Df^n(x)\overrightarrow{v}||_w &= ||Df(f^{n-1}(x))\cdot Df^{n-1}(x)\overrightarrow{v} ||_w\\
    &= ||Df(f^{n-1}(x))\cdot Df(f^{n-2}(x))\cdot Df^{n-2}(x)\overrightarrow{v} ||_w\\
    &= ||Df(f^{n-1}(x))\cdot Df(f^{n-2}(x))\cdots Df(x)\overrightarrow{v} ||_w
\end{align*}
Because $\mathcal{NW}(f)$ is forward invariant, we then have
\[||Df^n(x)\overrightarrow{v}||_w\geq \Lambda^n||\overrightarrow{v}||_w.\]
This implies
\begin{align}
    ||Df^n(x)\overrightarrow{v}||\geq \frac{\omega_0}{\omega_1} \Lambda^n||\overrightarrow{v}||.
\label{3}\end{align}

Because $TM_x=\mathbb{R}$ for all $x\in \mathcal{NW}(f)$, from $\eqref{3}$ we then get 
\[|Df^n(x)|\geq \frac{\omega_0}{\omega_1} \Lambda^n.\]

\begin{figure}[h]\label{fig1}
        \centering
  \includegraphics[width=0.41\linewidth]{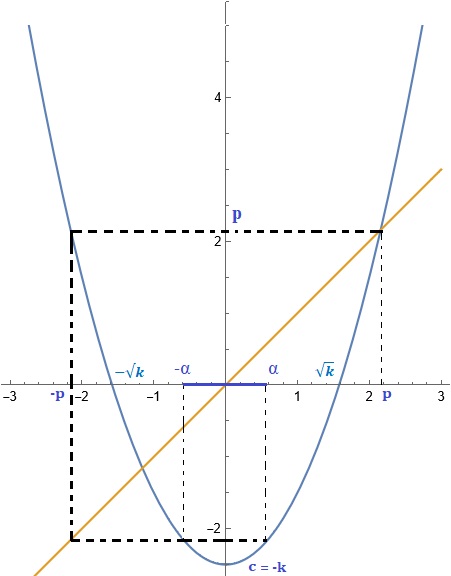}
       \caption{Graph of $f_c(x)=x^2+c$ when $c<-2$} \label{fig2}
            \end{figure}

Which proves the hyperbolicity of $f$. Thus if we can find a continuous function $\omega$ with the property $\eqref{3}$, then we are done. But notice that $\eqref{3}$ implies, 
\[2|x|\cdot \omega (f(x))\cdot ||\overrightarrow{v}|| \geq \Lambda \omega(x)||\overrightarrow{v}||. \] 
Therefore if we can find a continuous function $\omega$ such that 
\begin{align}
    2|x|\cdot \frac{\omega(f(x))}{\omega (x)}\gneq 1\ \text{for all}\ x\in \mathcal{NW}(f) \label{3.5}
\end{align}
then we are done.
For the convenience of calculations, let $k>0$ be such that $k=-c$ (Fig\ref{fig2}). Thus ${f(x)=x^2-k}$. Define the function ${G:\mathcal{NW}(f)\longrightarrow [0,\infty)}$ by
\[G(x):= \begin{cases}
       D((\alpha,p), (x,\sqrt{k})) &\quad\text{if }\ x\in \mathcal{NW}(f)\cap [\alpha, p],\\
       D((-p, -\alpha), (x,-\sqrt{k})) &\quad\text{if }\ x\in \mathcal{NW}(f)\cap [-p, -\alpha]. 
     \end{cases}\]
\newpage
\textbf{Remarks.}
\begin{enumerate}
    \item If $x=\pm \alpha, \pm p$, then we take $\limsup _{h\longrightarrow 0} G(x\pm h)$ accordingly.
    \item $\pm\sqrt{k}\notin \mathcal{NW}(f)$.
    \item $G(x)$ is continuous on $\mathcal{NW}(f)$.
\end{enumerate}

    \begin{figure}[h]\label{fig1}
        \centering
  \includegraphics[width=0.8\linewidth]{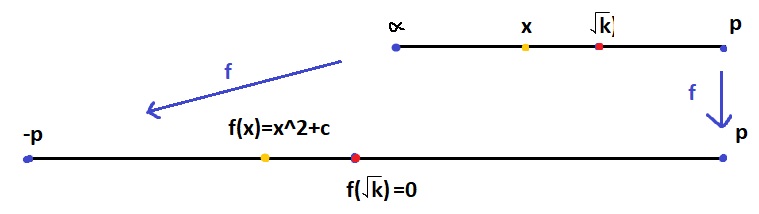}
        \caption{Graph of $f(x)=x^2-k$ maps $[\alpha,p]$ onto $[-p,p]$}
            \end{figure} 

Without the loss of generality let's assume $x\in \mathcal{NW}(f)\cap [\alpha, p]$. Then
\begin{align}
  \nonumber  \frac{G(f(x))}{G(x)} &=  \frac{D[f((\alpha,p)), f((x,\sqrt{k}))]}{D((\alpha,p), (x,\sqrt{k}))}\\
 \nonumber   &= \frac{2p|x^2-k|}{(x^2-k+p)p}\cdot \frac{(x-\alpha)(p-\sqrt{k})}{(p-\alpha)(|\sqrt{k}-x|)},\ \text{$\because\ p$ is a fixed point and $f(\sqrt{k})=0$.}\\
    &= \frac{2(p-\sqrt{k})}{(p-\alpha)}\cdot \frac{(x+\sqrt{k})(x-\alpha)}{x^2-k+p}
\end{align}

Let $g(x):=\frac{(x+\sqrt{k})(x-\alpha)}{x^2-k+p}$. Notice that $\alpha^2=k-p$. Thus $g(x)=\frac{(x+\sqrt{k})(x-\alpha)}{(x-\alpha)(x+\alpha)}$. 

First consider the case $x=\alpha$. By L'Hospital's rule,
\begin{align}
    g(\alpha)=\lim_{h\longrightarrow 0}\frac{(\alpha+h+\sqrt{k})h}{h(2\alpha +h)}= \frac{\alpha + \sqrt{k}}{2\alpha}.
\end{align}

Hence 
\begin{align}
    \nonumber \frac{G(f(\alpha))}{G(\alpha)}&=\frac{2(p-\sqrt{k})}{(p-\alpha)}\cdot \frac{\alpha + \sqrt{k}}{2\alpha}\\
    \nonumber &= \frac{(p-\sqrt{k})(\alpha+\sqrt{k})}{\alpha (p-\alpha)} \\
    \nonumber&= \frac{p\alpha-k+\sqrt{k}(p-\alpha)}{\alpha p - \alpha^2}\\
    \nonumber &= \frac{p\alpha-k+\sqrt{k}(p-\alpha)}{\alpha p - k + p},\ \because\ \alpha^2=k-p\\
    \nonumber &\geq 1 + \frac{\sqrt{k}(p-\frac{1}{2})-p}{\alpha p - k + p},\ \because\ \alpha \leq \frac{1}{2}\\
    \nonumber &\geq 1 + \frac{\sqrt{k}(2-\frac{1}{2})-p}{\alpha p - k + p}\\
    \nonumber &= 1+\frac{\frac{3}{2}\sqrt{k}-p}{\alpha p - k + p},\ \because\ p\geq 2\\
    &\gneq 1,\ \because p\geq 2\ \text{implies}\ \frac{3}{2}\sqrt{k}-p= \frac{3}{2}(p^2-p)-p > 0.\label{3.8}
\end{align}

Now suppose $x\gneq \alpha$. Then \[g(x)=\frac{(x+\sqrt{k})(x-\alpha)}{(x-\alpha)(x+\alpha)}= \frac{x+\sqrt{k}}{x+\alpha}\geq \frac{p+\sqrt{k}}{p+\alpha}.\]

Therefore form $(7)$,
\begin{align}
    \nonumber \frac{G(f(x))}{G(x)} &\geq \frac{2(p-\sqrt{k})}{p-\alpha} \cdot \frac{p+\sqrt{k}}{p+\alpha}\\
    \nonumber &= 2\cdot \frac{p^2-k}{p^2-\alpha^2}\\
    \nonumber &= 1,\ \because\ p=f(p)=p^2-k\ \text{and}\ \alpha^2=k-p\\
    &= \frac{G(f(p))}{G(p)}.\label{3.9}
\end{align}
 It is clear that $g$ is a decreasing function as $g'(x)<0$. Thus when when $x=p$ we get the minimum of $G(f(x))/G(x)$, which is $1$. Thus by $\eqref{3.8}$ and $\eqref{3.9}$ we have $G(f(x))/G(x)\geq 1$. The equality occurs if and only if $x=p$. So choose $M>0$, say, such that 
 \begin{align}
     \Big( \frac{G(f(\alpha))}{G(\alpha)}\Big)^M \gneq \frac{1}{2\alpha}.\label{3.10}
 \end{align} 
 Define
 \[\omega(x):=G(x)^M.\]
 Then by \eqref{3.8}, \eqref{3.9}, \eqref{3.10}, and the fact that $p\gneq 2$, we get
 
 \[2|x|\cdot \frac{ \omega(f(x))}{\omega(x)}\gneq 1 \]
 for any $x\in \mathcal{NW}(f)\cap [\alpha, p]$ \.\\ 
 
 Similarly we can extend $\omega$ to the whole $\mathcal{NW}(f)$ having the desired property. Therefore by \eqref{3.5}, we have the claim.

\ack 
I would like to thank my advisors Professor Bruce P. Kitchens and Professor Roland K. W. Roeder of Indiana University-Purdue University--Indianapolis (IUPUI) for their helpful comments and guidance.


\end{document}